

\input amstex.tex
\loadmsam
\loadmsbm
\loadbold
\input amssym.tex
\overfullrule=0pt
\baselineskip=13pt plus 2pt

\documentstyle{amsppt}
\pageheight{46pc}
\pagewidth{33pc}
\magnification=1200
\parskip 2pt
\NoBlackBoxes

\def\onto{\twoheadrightarrow}
\def\into{\hookrightarrow}
\def\s{\sigma}

\def\1{\bold 1}
\def\epsilon{\varepsilon}

\def\l{\langle} 

\def\hb{\hfil\break}
\def\ov #1{\overline{#1}}
\def\n{\noindent}
\def\bx{$\square$}

\def\MHS{\operatorname{MHS}}
\def\AHC{\operatorname{AHC}}

\def\br{\operatorname{bar}}

\def\Int{\operatorname{Int}}
\def\id{\operatorname{id}}

\def\Aut{\operatorname{Aut}}
\def\loc{\operatorname{loc}}

\def\Res{\operatorname{Res}}

\def\Gr{\operatorname{Gr}}
\def\vcd{\operatorname{vcd}}
\def\cd{\operatorname{cd}}
\def\char{\operatorname{char}}

\def\Gal{\operatorname{Gal}}

\def\Vec{\operatorname{Vec}}
\def\CM{\operatorname{CM}}
\def\ab{\operatorname{ab}}

\def\red{\operatorname{red}}
\def\Lie{\operatorname{Lie}}
\def\Inv{\operatorname{Inv}}

\def\Art{\operatorname{Art}}
\def\Hom{\operatorname{Hom}}

\def\DR{\operatorname{DR}}
\def\Sper{\operatorname{Sper}}
\def\Spec{\operatorname{Spec}}
\def\Arch{\operatorname{Arch}}
\def\ad{\operatorname{ad}}
\def\Rep{\operatorname{Rep}}
\def\Hod{\operatorname{Hod}}

\def\C{\Bbb C}
\def\G{\Bbb G}

\def\Q{\Bbb Q}
\def\R{\Bbb R}
\def\Z{\Bbb Z}

\loadeusm
\def\scr#1{{\fam\eusmfam\relax#1}}

\def\uM{\scr M}

\def\uMR{\scr M\scr R}
\def\l{\ell}

\def\unM{\underline M}

\def\unHom{Hom}

\def\tto{\tilde\to}

\def\u1{\bold 1}
\def\uHom{\Hom}
\def\uHod{\Hod}
\def\uAut{\Aut}
\def\uRep{\Rep}

\def\lim{\operatorname{lim}}
\def\liminv{\underset{\longleftarrow_N}\to\lim\,}
\def\limin{\underset{\longleftarrow_K}\to\lim\,}

\rightheadtext{Motivic torsors}

\topmatter
\title Motivic torsors\endtitle
\author Yuval Z. FLICKER\endauthor
\footnote"~"{Department of Mathematics, The Ohio State University,
231 W. 18th Ave., Columbus, OH  43210-1174;\hb
flicker\@math.ohio-state.edu; \, 
1991 Mathematics Subject Classification: 14L30; 14F20, 14F40, 14P25, 11G10.}

\abstract The torsor $P_\s=\uHom^\otimes(H_{\DR},H_\s)$ under the motivic 
Galois group $G_\s=\uAut^\otimes H_\s$ of the Tannakian category $\uM_k$ 
generated by one-motives related by absolute Hodge cycles 
over a field $k$ with an embedding $\s:k\into\C$ is shown to be determined 
by its projection $P_\s\to P_\s/G_\s^0$ to a $\Gal(\ov k/k)$-torsor, and by 
its localizations $P_\s\times_kk_\xi$ at a dense subset of orderings $\xi$ 
of the field $k$, provided $k$ has virtual cohomological dimension (vcd) one.
This result is an application of a recent local-global principle for
connected linear algebraic groups over a field $k$ of vcd $\le 1$. 
\endabstract
\endtopmatter
\document

The singular cohomology with coefficients in the field $\Q$ of rational
numbers of a smooth projective -- even just complete -- variety over $\C$ has 
a (``pure'') Hodge structure. Motives with a realization (usually by means of 
some cohomology theory) which has a pure Hodge 
structure are called pure motives. Deligne defined in [D-II] a mixed Hodge 
structure to be a finite dimensional vector space $V$ over $\Q$ with a finite 
increasing (weight) filtration $W_{\bullet}$ and a finite decreasing 
(Hodge) filtration $F^{\bullet}$ on $V\otimes_\Q\C$ such that $F^{\bullet}$ 
induces a Hodge structure of weight $n$ on the graded piece $\Gr_n^WV
=W_nV/W_{n-1}V$ for each $n$. Deligne showed in [D-III] that the cohomology 
$H^\ast(E(\C),\Q)$ of any variety $E$ over $\C$ -- not necessarily complete 
and smooth -- carries a natural mixed Hodge structure. Motives with a 
realization which has a mixed Hodge structure are called mixed motives 
for emphasize.

Deligne introduced the notion of a one-motive $M$ -- as well as its
dual $M^{\sssize\vee}$, and Betti: $M(\C)_B$, de Rham: $H_{\DR}(M)$, and 
$\l$-adic: $H_{\l}(M)$, realizations -- in [D-III], \S 10, as a simple 
example of a motive whose Betti realization $M(\C)_B$ has a mixed Hodge 
structure, but does not have a Hodge structure. 
Let $\s:k\into\C$ be an embedding of a field $k$ in the field $\C$ of 
complex numbers, and $\ov\s:\ov k\into\C$ an extension to an algebraic
closure $\ov k$. Write $\Gal(\ov k/k)$ for the Galois group. For a variety
$E$ over $k$, write $\s E$ for the $\C$-variety $E\times_{k,\s}\C$.

A one-motive over $k$ is a complex $M=[X\overset u\to\to E]$ of length one 
placed in degrees $0$ and $1$, comprising of a semi-abelian variety $E$ 
(namely an extension $1\to T\to E\to A\to 0$ of an abelian variety $A$ by a 
torus $T$) over $k$, a finitely generated torsion free $\Gal(\ov k/k)$-module 
$X$, and a $\Gal(\ov k/k)$-equivariant homomorphism $u:X\to E(\ov k)$. Note
that $E$ is a commutative $k$-group.
One-motives include the Artin motives as $[X\to 0]$ and the Tate motive as 
$[0\to\G_m]$. We also write $M=(X,A,T,E,u)$, $M\otimes\Q$ for the isogeny 
class of $M$, $\s M=[X\overset u\to\to \s E]$ and 
$\s M(\C)=[X\overset u\to\to \s E(\C)]$. A one-motive $M$ has a 
``weight'' filtration: $W_0M=[X\overset u\to\to E]$, $W_{-1}M=[0\to E]$, 
$W_{-2}M=[0\to T]$, $W_{-3}M=[0\to 0]$, with graded factors $\Gr_0M=X$, 
$\Gr_{-1}M=(E/T)[-1]=A[-1]$, and $\Gr_{-2}M=T[-1]$. 
Put $\Gr^WM=[X\overset 0\to\to A\times_k T]$.

The Betti realization $H_{\s}(M)=\s M(\C)_B$ of a one-motive $M=[X\overset 
u\to\to E]$ over $k$ is the vector space $T_{\s}(M)\otimes\Q$, where the 
lattice $T_{\s}(M)$ is the fiber product of $\Lie \s E(\C)$ and $X$ over 
$\s E(\C)$, namely the pullback of $0\to H_1(\s E(\C))$ $\to\Lie \s E(\C)$
$\overset\exp\to\to \s E(\C)\to 1$ by $X\overset u\to\to \s E(\C)$. It depends
on the embedding $\s:k\into\C$. Then $\s M(\C)_B$ is a mixed Hodge structure 
$(V,W_{\bullet},F^{\bullet})$ of type $\{(0,0),\,(0,-1),\,(-1,0),\,(-1,-1)\}$ 
whose graded parts are $\Gr_0^WV=X\otimes\Q$, polarizable $\Gr_{-1}^WV=
H_1(\s A(\C),\Q)$, and $\Gr_{-2}^WV=H_1(\s T(\C),\Q)=X_\ast(\s T)\otimes\Q$; 
see [D-III], 10.1.3. 

Denote by $\uM_k$ the Tannakian category (Deligne-Milne [DM], Definition 2.19)
generated by the isogeny classes of one-motives over $k$, in the category
$\uMR_k$ of mixed realizations (Jannsen [J], 2.1), related by absolute Hodge
cycles (Deligne [D2], 2.10, Brylinski [Br], 2.2.5). The objects of $\uMR_k$
are tuples $H=(H_{\DR},H_\l,H_\s;I_{\infty,\s},I_{\l,\ov\s})$, where $\l$
ranges over the rational primes, $\s$ over the embeddings $k\into\C$, and
$\ov\s$ over the $\ov k\into\C$, described in [J], p. 10. In particular 
$H_{\DR}$ is a finite dimensional $k$-vector space with a decreasing (Hodge)
filtration $(F^n;n\in\Z)$ and an increasing (weight) filtration 
$(W_m;m\in\Z)$; $H_\l$ is a finite dimensional $\Gal(\ov k/k)$-module over 
$\Q_l$ with $\Gal(\ov k/k)$-equivariant increasing (weight) filtration 
$W_{\bullet}$; $H_\s$ is a mixed Hodge structure (over $\Q$), and 
$I_{\infty,\s}:H_\s\otimes\C\tto H_{\DR}\otimes_{k,\s}\C$, 
$I_{\l,\ov\s}:H_\s\otimes\Q_\l\tto H_{\l}$ ($\s=\ov\s\vert k$) 
are the comparison isomorphisms. 

The morphisms in $\uMR_k$ are tuples $(f_{\DR},f_\l,f_\s)_{\l,\s}$ 
described in [J], p. 11. In particular $f_\s:H_\s\to H'_\s$ is a 
morphism of mixed Hodge structures, $f_{\DR}:H_{\DR}\to H'_{\DR}$ is 
$k$-linear and $f_\l:H_\l\to H'_\l$ is a $\Q_\l$-linear
$\Gal(\ov k/k)$-morphism, which correspond under the comparison isomorphisms.
The category $\uMR_k$ is abelian ([J], 2.3), tensor ([J], 2.7) with identity
object $\1=(k,\Q_\l,\Q;\id_{\infty,\s},\id_{\l,\ov\s})$, and it has internal 
$\unHom(H,H')\in\uMR_k$ for all $H,\,H'$ in $\uMR_k$ (thus 
$\Hom(H'',\unHom(H,H'))=\Hom(H''\otimes H,H')$ for all 
$H,\,H',\,H''\in\uMR_k$). For example, $H_{\DR}(\unHom(H,H'))
=\Hom_k(H_{\DR},H'_{\DR})$, $H_{\l}(\unHom(H,H'))=\Hom_{\Q_\l}(H_\l,H'_\l)$, 
$H_{\s}(\unHom(H,H'))=\Hom_{\Q}(H_\s,H'_\s)$. Hence $\uMR_k$ is rigid (each 
object $H$ has a dual $H^{\sssize\vee}=\unHom(H,1))$.

Defining the space $\AHC(H)$ of absolute Hodge cycles of $H\in\uMR_k$ to be 
the set of $(x_{\DR},x_\l,x_\s)\in H_{\DR}\times\prod_\l H_\l\times\prod_\s 
H_\s$ such that $I_{\infty,\s}(x_\s)=x_{\DR}$, $I_{\l,\ov\s}(x_\s)=x_\l$
for all $\s$, $\ov\s$ with $\ov\s|k=\s$ and $x_{\DR}\in F^0H_{\DR}\cap W_0
H_{\DR}$ (it is a finite dimensional vector space over $\Q$), one has 
$\Hom(H,H')=\AHC(\unHom(H,H'))$. A Hodge cycle with respect to $\s$ is a tuple
$(x_{\DR},x_\l)\in H_{\DR}\times\prod_\l H_\l$ such that there is $x_\s\in 
H_\s$ with $I_{\infty,\s}(x_\s)=x_{\DR}$, $I_{\l,\ov\s}(x_\s)=x_\l$, $x_{\DR}
\in F^0H_{\DR}\cap W_0 H_{\DR}$. Then $\uMR_k$ is a Tannakian category neutral
over $\Q$, namely a rigid abelian tensor $\Q$-linear category with a 
$\Q$-valued fiber ([DM], Definition 2.19: exact faithful $\Q$-linear tensor) 
functors $H^{\#}_\s:\uMR_k\to\Vec_{\Q}$, $H\mapsto H^{\#}_\s$. The $\#$ 
emphasizes here that the symbol indicates the underlying vector space.
In the literature, and in the abstract of this paper, $\#$ is omitted to
simplify the notations for the reader who knows when $H_\s$ is regarded as
a mixed Hodge structure, and when it is regarded only as a vector space.

The mixed realization $H(M)$ of a one-motive $M$ is 
$(H_{\DR}(M),H_\l(M),H_\s(M);\,I_{\infty,\s},I_{\l,\ov\s})$; see [D-III],
10.1.3: the $H$ are $H_1$. Note that the dual one motive $M^{\sssize\vee}$ 
(introduced in [D-III], 10.2.11) satisfies 
$H(M^{\sssize\vee})=\unHom(H(M),\Q(1))$.
Hence $H(M)^{\sssize\vee}=H(M^{\sssize\vee})(-1)$. From now on by a motive
we mean an object in the Tannakian category $\uM_k$ generated in $\uMR_k$ by 
one-motives. The functor $H^{\#}_\s$ -- which associates to a motive $M$ the 
vector space underlying the mixed Hodge Betti realization $\s M(\C)_B$ -- is 
a fiber functor on $\uM_k$, making $\uM_k$ Tannakian and neutral over $\Q$. 
Note that an isomorphic -- but not canonically -- fiber functor is $H^{\#}_\s
\Gr^W$. This fiber functor corresponds to a choice of a Levi decomposition of 
the motivic Galois group, see the end of the 5th paragraph below.

The category $\uM_k$ is not semi-simple, but it has a semi-simple Tannakian 
full subcategory $\uM_k^{\red}$ of motives generated by abelian varieties 
($M=[0\to A]$) and Artin motives ($M=[X\to 0]$) over $k$, related by absolute 
Hodge cycles ([DM], Propositions 6.5 and 6.21). Thus it is the subcategory of 
$\uMR_k$ generated by $H(A)(=(H_{1,\DR}(A),H_{1,\text{\'et}}(A\times_k\ov k,
\Q_{\ell}),$ $H_1(\s A(\C),\Q))$ of the abelian varieties $A$ over $k$, and 
the Artin motives $H(X)=X\otimes\1=(X\otimes k,X\otimes\Q_{\ell},X\otimes\Q)$.
Note that the realization $H(T)$ of the torus $[0\to T]$ is the Tate twisted 
Artin motive $X_\ast(T)\otimes\1(1)
(=(X_\ast(T)\otimes k(1),X_\ast(T)\otimes
\Q_{\ell}(1),X_\ast(T)\otimes\Q(1)))$,
where $X_\ast(T)=\unHom(\G_m,T)$ (internal $\Hom$ in the category of
one-motives). The subcategory $\uM_k^{\red}$ 
of $\uM_k$ is also neutral over $\Q$, by the fiber functor $H^{\#}_\s$. 

Denote by $\uM_k\otimes k$ the category $(\uM_k)_{(k)}$ of [DM], Proposition 
3.11, obtained on extending coefficients from $\Q$ to $k$. It is a Tannakian 
category neutral over $k$. The functors $H^{\#}_\s\otimes k$ 
$(:M\mapsto \s M(\C)_B\otimes k)$ and $H^{\#}_{\DR}$ on $\uM_k\otimes k$ 
are fiber functors with values in $k$. The groups 
$G_\s=\uAut^\otimes (H^{\#}_\s\otimes k|\uM_k\otimes k)$ and 
$G_{\DR}=\uAut^\otimes(H^{\#}_{\DR}|\uM_k\otimes k)$ 
of automorphisms of the fiber functors are affine group schemes over $k$ 
([DM], Theorem 2.11 and Proposition 3.11); they are inner forms of each other.
Even a conjectural description of these groups is elusive. The functors 
$H^{\#}_\s\otimes k$ and $H^{\#}_{\DR}$ define equivalences 
$\uM_k\otimes k\tto$ $\uRep_kG_\s$ and $\uM_k\otimes k\tto\uRep_kG_{\DR}$ 
of tensor categories. 

Similarly we have the Tannakian category $\uM_k^{\red}\otimes k$, 
which is semi-simple and neutral over $k$ by the fiber functors
$H^{\#}_\s\otimes k$ and $H^{\#}_{\DR}$, the $k$-groups 
$G_\s^{\red}=\uAut^\otimes (H^{\#}_\s\otimes k|\uM_k^{\red}\otimes~k)$ 
and $G_{\DR}^{\red}=\uAut^\otimes(H^{\#}_{\DR}|\uM_k^{\red}\otimes k)$, 
and the equivalences 
$\uM_k^{\red}\otimes k\tto$ $\uRep_kG_\s^{\red}$ and 
$\uM_k^{\red}\otimes k\tto\uRep_kG_{\DR}^{\red}$. Since the category 
$\uM_k^{\red}\otimes k$ is semi-simple, [DM], Remark 2.28 implies that
$G_\s^{\red}$ and $G_{\DR}^{\red}$ are pro-reductive (meaning that the
connected component is the projective limit of connected reductive groups).
The group $G_\s^{\red}$ (resp. $G_{\DR}^{\red}$) is the maximal pro-reductive 
quotient of the affine group scheme $G_\s$ (resp. $G_{\DR}$). 

Note that a $\otimes$-functor $F:A\to B$ of Tannakian categories
and a fiber functor $\beta$ on $B$ define a map $f:G_B=\Aut^\otimes(\beta)
\to G_A=\Aut^\otimes(\beta\circ F)$ of the motivic groups (the image 
$g^A=(g^A_{X_A})=f(g^B)$ is defined by $g^A_{X_A}=g^B_{F(X_A)}$), and vice 
versa: $f:G_B\to G_A$ defines $F:A=\Rep G_A\to B$. For relations of 
properties of $F$ and $f$ see Saavedra [Sa], II, 4.3.2. 

Denote by $U_\s$ the kernel of the projection $G_\s\to G_\s^{\red}$; it is 
pro-unipotent. By the Levi decomposition, the extension 
$1\to U_\s\to G_\s\to G_\s^{\red}\to 1$ 
splits. More precisely, the essentially surjective functor (a functor is 
called {\it essentially surjective} if each object in the target category 
is isomorphic to an object in the image of the functor) 
$\Gr^W:\uM_k\to\uM_k^{\red}$, defined on one-motives by 
$M=(X,A,T,E,u)\mapsto H(X)\oplus H(A)\oplus H(X_\ast(T))(1)$, is an inverse 
to $\uM_k^{\red}\into\uM_k$. Correspondingly $G_\s^{\red}=\uAut^\otimes 
(H^{\#}_\s\otimes k|\uM_k^{\red}\otimes~k)$ is canonically a subgroup of
$\Gr^WG_\s=\Aut^\otimes(H^{\#}_\s\Gr^W\otimes k|\uM_k\otimes k)$, which is 
isomorphic by the Levi decomposition -- but not canonically -- to 
$G_\s=\Aut^\otimes(H^{\#}_\s\otimes k|\uM_k\otimes k)$.

Our main object of study is the affine scheme 
$P_\s=\uHom^\otimes(H^{\#}_{\DR},
H^{\#}_\s\otimes k;\uM_k\otimes k)$ over $k$ of morphisms of fiber functors 
([DM], Theorem 3.2). It is a $G_\s$-torsor (right principal homogeneous space)
over $k$, and so it defines a class $h_\s$ of the first Galois cohomology set 
$H^1(k,G_\s)=H^1(\Gal(\ov k/k),G_\s(\ov k))$. The group $G_\s$ is called the 
($\s$-){\it motivic Galois group} of $\uM_k\otimes k$, and $P_\s$ the 
($\s$-){\it motivic torsor} of $\uM_k\otimes k$. Analogously we have the 
$G_\s^{\red}$-torsor $P_\s^{\red}=\uHom^\otimes(H^{\#}_{\DR},H^{\#}_\s
\otimes k;\uM_k^{\red}\otimes k)$ 
over $k$, and its class $h_\s^{\red}$ in $H^1(k,G_\s^{\red})$. The 
$G_\s^{\red}$-torsor $P_\s^{\red}$ is the quotient $P_\s/U_\s$. 

Denote by $\uM_k^0$ the Tannakian subcategory generated by Artin motives 
$[X\to 0]$ in 
$\uM_k$. It is equivalent to the category of [DM], Proposition 6.17, 
generated by the zero dimensional varieties $Z$ over $k$. The motivic Galois 
group $\uAut^\otimes(H^{\#}_\s\otimes k|$ $\uM_k^0\otimes k)$ of 
$\uM_k^0\otimes k$ is the constant pro-finite group scheme 
$\Gamma_k=\limin\coprod_\gamma(\Spec k)_\gamma$ [($k\subset) K$ finite Galois 
extensions, $\gamma\in\Gal(K/k)$] over $k$ (with structure morphisms 
$\coprod_{\gamma\in\Gal(K/k)}\id_\gamma$). Its group of $\ov k$-points is 
$\Gal(\ov k/k)$, and the functor $H^{\#}_\s\otimes k$ ($:X\mapsto X\otimes k$,
or $:Z\mapsto k^{Z(\ov k)}$ in [DM],
6.17) induces an isomorphism $\uM_k^0\otimes k$ 
$\tto\uRep_k(\Gamma_k)$ ([DM], Proposition 6.17). The group $\Gamma_k$ 
is the group of connected components of $G_\s^{\red}$ ([DM], Proposition 
6.23(a,b)). [Note that the proofs of Propositions 6.22(a), 6.23 of [DM] are 
incorrect for the full category of pure motives as stated there, but they do 
apply in our context of motives of abelian varieties and one-motives; see 
Remark 1 at the end of this paper.] 

Thus the inclusion $\uM_k^0\into\uM_k^{\red}$ defines a surjection 
$G_\s^{\red}\overset\pi\to\onto\Gamma_k$ (by [DM], Remark 2.29). Its kernel 
$G_\s^{\red,0}$ is the connected component of the identity of $G_\s^{\red}$, 
a connected pro-reductive affine $k$-group scheme which is the motivic Galois 
group $\uAut^\otimes(H^{\#}_\s\otimes k|\uM_{\ov k}^{\red}\otimes k)$ of 
$H^{\#}_\s\otimes k$ on $\uM_{\ov k}^{\red}\otimes k$. The almost surjective 
functor (we say that a functor is {\it almost surjective} if each object of 
the target category is isomorphic to a subquotiet of an object in the image 
of the functor; see [DM], Proposition 2.21(b))  
$\uM_k^{\red}\to\uM_{\ov k}^{\red}$, $A\mapsto\ov A=A\times_k\ov k$, defines 
the injection $G_\s^{\red,0}\overset\iota\to\to G_\s^{\red}$. In particular, 
denote the quotient $p:P_\s^{\red}\to P_\s^{\red}/G_\s^{\red,0}$ by 
$P_\s^{\Art}$. It is the $\Gamma_k$-torsor $\uHom^\otimes(H^{\#}_{\DR},
H^{\#}_\s\otimes k;\uM_k^0\otimes k)$. Its class $h_\s^{\Art}$ in 
$H^1(k,\Gamma_k)=H^1(\Gal(\ov k/k),\Gamma_k(\ov k))$ is the image of 
$h_\s=\{P_\s^{\red}\}$ under the map $H^1(k,G_\s^{\red})\to H^1(k,\Gamma_k)$.

Since $G_\s$ is the semi-direct product of the pro-reductive $G_\s^{\red}$
and the pro-unipotent $U_\s$, we have that $\Gamma_k$ is the group of 
connected components of $G_\s$. The inclusion $\uM_k^0\to\uM_k$ defines a 
surjection $G_\s\overset\pi\to\to\Gamma_k$ ([DM], Proposition 2.21(a)), whose 
kernel $G_\s^0$ is the connected component of the identity of $G_\s$. This 
connected affine $k$-group scheme is the motivic Galois group 
$\uAut^\otimes(H^{\#}_\s\otimes k|\uM_{\ov k}\otimes k)$ of 
$H^{\#}_\s\otimes k$ on $\uM_{\ov k}\otimes k$. The almost surjective 
functor $\uM_k\to\uM_{\ov k}$, defined on one-motives by 
$M\mapsto\ov M=M\times_k\ov k=[X\overset u\to\to E\times_k\ov k]$, 
induces the injection $G_\s^{0}\overset\iota\to\to G_\s$. The quotient 
$p:P_\s\to P_\s/G_\s^{0}$ is $P_\s^{\Art}$. Its class in $H^1(k,\Gamma_k)$ is 
the image of $h_\s=\{P_\s\}$ under the map $H^1(k,G_\s)\to H^1(k,\Gamma_k)$.
The functor $H_{\s}\otimes k$ maps the Tannakian category $\uM_{\ov k}\otimes 
k$ to the Tannakian category of $k$-mixed Hodge structures. This would help
us understand what we need to know about our motivic objects, but this map 
is not fully faithful when $k\not=\Bbb Q$.

The statement of our theorem uses the set $\Sper k$ of orderings $\xi$ of the 
field $k$. It is a compact totally disconnected topological space, where a 
basis of the topology is given by the sets $\{\xi;\,a>0$ in $\xi\}$ for all 
$a$ in $k$ (see, e.g., Scharlau [Sc], Ch. 3, \S 5). The space $\Sper k$ is 
naturally homeomorphic to the quotient of the space $\Inv(\Gal(\ov k/k))$ of 
involutions (elements of order precisely two) in $\Gal(\ov k/k)$ (endowed 
with the usual profinite topology) by conjugation under $\Gal(\ov k/k)$. 
Denote by $k_\xi$ a real closure of $k$ (in $\ov k\subset\C$) whose ordering 
induces $\xi$ on $k$. Then $\Gal(\ov k/k_\xi)$ is generated by $c_\xi$ in 
$\Inv(\Gal(\ov k/k))$. If $c$ is an involution in $\Gal(\ov k/k)$, for any 
field $k$, then $\char k=0$, the fixed field of $c$ in $\ov k$ is a real 
closure $k_\xi$ of $k$ whose ordering induces $\xi$ on $k$, 
$\ov k=k_\xi(\sqrt{-1})$, and the restriction of $c$ to the algebraic closure 
$\ov\Q$ of $\Q$ is non trivial (it  takes $\sqrt{-1}$ to $-\sqrt{-1}$),
hence it is in the unique conjugacy class of involutions in $\Gal(\ov\Q/\Q)$.
The ordered field $k$ (or $k_\xi$) embeds in a real closed field $R_\xi$ of
codimension $2$ in $\C$ -- thus $\C=R_\xi(\sqrt{-1})$ -- whose ordering
induces $\xi$ on $k$. An ordering $\xi$ of $k$ is called {\it archimedean} if 
the real closure $k_\xi$ embeds in $\R$. When $k$ is finitely generated, the 
set $\Arch k$ of archimedean orderings in $k$ is dense in $\Sper k$; this is
shown below.

The affine $k_\xi$-scheme $P_{\s,\xi}=P_\s\times_kk_\xi$ is a $G_\s$-torsor 
over $k_\xi$. Its class $h_{\s,\xi}$ in 
$H^1(k_\xi,G_\s)=H^1(\Gal(\ov k/k_\xi),
G_\s(\ov k))$ is the image of $h_\s$ under the natural localization map 
$H^1(k,G_\s)$ $\to$ $ H^1(k_\xi,G_\s)$. 
Alternatively it can be described using
the fact that the natural map $H^1(k_\xi,G_\s)$ 
$\to$ $ H^1(R_\xi,G_\s)$ is an 
isomorphism (this is implied by the Artin-Lang theorem 
(see [BCR], Th\'eor\`eme
4.1.2)), as follows. The continuous map 
$\s M(\C)\to \s M(\C)$ ($M\in\uM_k$) 
defined by $c_\xi\not=1$ in $\Gal(\C/R_\xi)$ 
induces an involutive endomorphism
of $\s M(\C)_B$. The image in $G_\s(\C)$ defines a (Galois) cohomology class
in $H^1(R_\xi,G_\s)$, which is $h_{\s,\xi}$.

Let $k$ be a field with virtual cohomological dimension $\le 1$ (thus 
$\vcd(k)=\cd(k(\sqrt -1))$ is at most one). We have $\vcd(k)=\cd(k)$ precisely
when $k$ has no orderings, thus $\Sper k$ is empty. Examples of $k$ with
$\vcd k=1<\cd k$ are $k=\R(x)$ or $R(x)$, where $R$ is a real closed
field (Serre [S1], II, 3.3(b)), $\R((x))$ and $R((x))$ ([S1], II, 3.3, Ex. 
3), and $\Q^{\ab}\cap\R$ ([S1], II, 3.3, Proposition 9). We assume that $k$
embeds in $\C$ (to use [DM]; to embed in $\C$ a field $k$ of 
cardinality bounded by that of $\C$, choose transcendence bases in both).
Fix $\s:k\into\C$.

\proclaim{Theorem} Let $p':P'\to P_\s^{\Art}$ be a $G_\s$-torsor over $k$ such
that $P'_\xi=P'\times_kk_\xi$ is isomorphic to $P_{\s,\xi}=P_\s\times_kk_\xi$ 
for $\xi$ in a dense subset of $\Sper k$. Then there exists an isomorphism of 
$G_\s$-torsors $\lambda:P_\s\to P'$ such that $p'\circ\lambda=p$.

The same result holds with $G_\s$ and $P_\s$ replaced by $G_\s^{\red}$ and 
$P_\s^{\red}$.
\endproclaim

Our work is influenced by Blasius-Borovoi [BB] who considered the number
field $\Q$ (whose $\vcd$ is $2$) and the semi-simple Tannakian subcategory 
$\uM^{\red,H}_\Q$ generated by Artin motives and motives of abelian varieties
$A$ over $\Q$ for which the 
group $((G_\s^A)^0)_\Bbb R^{\ad}$ has no factor of type $D_n^H$ (in the
notations of Deligne [D1], (1.3.9)), and by 
Wintenberger [W] who had considered
the field $\Q$ and the semi-simple Tannakian subcategory $\uM_\Q^{\red,\CM}$ 
generated by Artin motives and motives of abelian varieties with complex 
multiplication over $\Q$. 

Our theorem is an application of the local-global principle for a field $k$ 
with $\vcd(k)\le 1$. We can work in the generality of the entire category 
$\uM_k$ and the group $G_\s$ by virtue of the local-global principle:
$H^1(k,G)\hookrightarrow\prod_\xi H^1(k_\xi,G)$, proven by Scheiderer [Sch] 
for a perfect field $k$ with $\vcd k\le 1$ and a connected $k$-linear 
algebraic group $G$. In the number field case the analogous well known 
local-global principle holds only for semi-simple simply connected $G$. 

When $\cd(k)\le 1$, 
thus when $k$ has no orderings, the class of $P_\s$ is determined by 
$P_\s^{\Art}$ alone. To deal with this case, we use only Steinberg's theorem 
([S1], III, \S 2.3) on the vanishing of $H^1(k,G)$ for a perfect field $k$ 
with $\cd(k)\le 1$ and a connected $k$-linear algebraic group $G$. 

It will be interesting to study our motivic objects over fields $k$ with $\vcd
\le 2$. In this context, note that a local-global principle for $k$ with 
$\vcd(k)\le 2$ and semi-simple simply connected classical linear algebraic 
groups has recently been established by Bayer-Fluckiger and Parimala [BP].

It is my pleasure to express my deep gratitude to P. Deligne for watching over
my first steps in the motivic fairyland, to M. Borovoi, U. Jannsen, R. Pink,
C. Scheiderer, J.-P. Serre, R. Sujatha, and the Referee, for useful comments, 
to M. Jarden for invitation to talk on this work at the Gentner Symposium 
on Field Arithmetic, Tel-Aviv University, October 1997, and to the National
University of Singapore for its hospitality in late 1999 while this paper was 
refereed. NATO grant CRG 970133 is gratefully acknowledged.

\demo{Proof of theorem} It is easy to adapt the proof to the context of the
pro-reductive quotient group $G_\s^{\red}$, so we discuss only the general
case of the entire group $G_\s$. 

Let $z\in Z^1(k,G_\s)$ be a $1$-cocycle representing $h_\s=\{P_\s\}\in 
H^1(k,G_\s)$. As in [S1], I.5.3, denote by ${}_zG_\s$ the form of 
$G_\s$ twisted by $z$. It is the affine group scheme over $k$ on which
$\Gal(\ov{k}/k)$ acts by $s:g\mapsto(\Int(z_s))(s(g))$ ($g\in 
G_\s(\ov{k})$, $s\in\Gal(\ov{k}/k)$). The natural bijection 
$H^1(k,{}_zG_\s)\tto H^1(k,G_\s)$, defined by $(x_s)\mapsto(x_sz_s)$ ([S1], 
I.5.3, Proposition 35), takes the trivial element of $H^1(k,{}_zG_\s)$ to 
$h_\s$. Denote by $\eta$ the class in $H^1(k,{}_zG_\s)$ which maps to $h'$, 
the class in $H^1(k,G_\s)$ of the $G_\s$-torsor $P'$. By the very definition 
of $P_\s$, as relating $G_\s$ and $G_{\DR}$, we have that $G_{\DR}$ is 
${}_{P_\s}G_{\s}=P_{\s}\times_{G_\s}G_\s$ (this is ${}_FP=P\times^AF$ in 
the notations of the first paragraph of [S1], I, \S 5.3; here $A$ of [S1] 
is $G_\s(\ov{k})$, which acts on $P_\s(\ov{k})$ (=$P$ in [S1])
by right multiplication and on $G_\s(\ov{k})$ (=$F$ in [S1]) by conjugation).
By the third paragraph of [S1], I, \S 5.3, we have that $G_{\DR}$ is the twist
${}_zG_\s$ of $G_\s$ by $z$. Since $P'_\xi\simeq P_{\s,\xi}$, the localization
$\eta_\xi=\loc_\xi(\eta)$ of $\eta$ in $H^1(k_\xi,G_{\DR})$ is $1$, for a 
dense set of $\xi$ in $\Sper k$. Since $P_\s$ and $P'$ project to the same 
$\Gamma_k$-torsor $P_\s^{\Art}$ in $H^1(k,\Gamma_k)$, the image of $\eta$ in 
$H^1(k,{}_{z'}\Gamma_k)$ is $1$, where $z'$ in $Z^1(k,\Gamma_k)$ is the image 
of $z\in Z^1(k,G_\s)$ under the projection $G_\s\to\Gamma_k$. Our aim is to 
show that $\eta=1$ in $H^1(k,G_{\DR})$.

Consider the exact sequence of affine group schemes
$$1\to G_{\DR}^0={}_zG_\s^0\to G_{\DR}={}_zG_\s\to\Gamma_{k,{\DR}}
={}_{z'}\Gamma_k\to 1.$$
Since the image of $\eta\in H^1(k,G_{\DR})$ in $H^1(k,\Gamma_{k,{\DR}})$ is
trivial, there is $\eta^0\in H^1(k,G_{\DR}^0)$ which maps to $\eta$. 
The group $G_{\DR}^0$ is a connected pro-finite affine group scheme
over $k$ ([DM], Proposition 6.22(a)). Thus $G_{\DR}^0=\liminv(G_{\DR}^N)^0$,
where $G_{\DR}^N$ is the motivic Galois group $\uAut^\otimes(H^{\#}_{\DR}|$
$\uM_{k_N}^N\otimes k_N)$ of the Tannakian subcategory $\uM_{k_N}^N$ of 
$\uM_k$  generated by a finite set $N$ of one-motives and their duals,
the Artin motives and the Tate motive $T$ and its dual $T^{\sssize\vee}$.
The finite set $N$ is defined over a finitely generated over $\Q$ subfield
$k_N$ of $k$. 

As explained in the proof of [DM], Proposition 6.22(a), 
$(G_{\DR}^N)^0$ is a linear algebraic group. Correspondingly $\eta^0=\liminv
\eta_N^0$, where $\eta_N^0\in H^1(k,(G_{\DR}^N)^0)$. Further, $\eta=\liminv
\eta_N$, where $\eta_N$ is the image of $\eta_N^0$ under the map
$H^1(k,(G_{\DR}^N)^0)\to H^1(k,G_{\DR}^N)$. Since $\eta_\xi$ is trivial in
$H^1(k_\xi,G_{\DR})$, the localization $\eta_{N,\xi}=\loc_\xi(\eta_N)$ is
trivial for all $N$, for the dense set of $\xi$ in $\Sper k$ of the theorem.

Write $\Arch k$ for the set of archimedean orderings in $\Sper k$.
The proposition below asserts that the homomorphism $G_{\DR}(k_\xi)\to
\Gamma_{k,{\DR}}(k_\xi)$ is surjective for every $\xi\in\Arch k$. In
particular $G_{\DR}^N(k_\xi)\onto\Gamma_{k_N,{\DR}}(k_\xi)=\Z/2$ for each
finite $N$ and $\xi\in\Arch k$. We claim that this map is onto for all
$\xi\in\Sper k$. The $k_N$-group $G^N_{\DR}$ has two connected components;
denote by $C=G^{N,+}_{\DR}$ the component not containing the identity. The
surjectivity means that $C(k_\xi)$ is non empty (for all 
$\xi\in\Arch k$). It follows from the Artin-Lang theorem that $C(k_{N,\xi})$
is non empty for all $\xi\in\Arch k_N$. But the set of $\xi\in\Sper k_N$
such that $C(k_{N,\xi})$ is non empty is open and closed in $\Sper k_N$ (see,
e.g., [Sch], Corollary 2.2). Our claim follows once we show that for a
finitely generated field $k_N$, the set $\Arch k_N$ is dense in $\Sper k_N$.

To show that for a finitely generated field $k_N$ the set $\Arch k_N$ is 
dense in $\Sper k_N$, choose a purely transcendental extension 
$F=\Q(t_1,\dots,t_n)$ of $\Q$ of finite codimension in $k_N$. Since the 
restriction of orderings is an open map $\Sper k_N\to\Sper F$, and an
ordering of $k_N$ is archimedean if its restriction to $F$ is, it suffices
to show that $\Arch F$ is dense in $\Sper F$. For this, we proceed to show
that the non empty basic open set defined by $p_1,\dots,p_r\in F$ contains
an archimedean ordering. The open set being non empty means that there is an 
ordering of $F$ which makes the $p_j$ positive.
In other words, there are a real closed field $R$ and $x\in R^n$
such that $p_j(x)>0$, all $j$. Then the same is true for $R=\R$,
by the Tarski principle (see, e.g., [BCR], I.1.4). That is, there is
$x\in\R^n$ such that $p_j(x)>0$, all $j$. The inequalities remain true in a
neighborhood of $x$, hence the components $x_1,\ldots,x_n$ of $x$ can be
chosen to be algebraically independent. The embedding $F\into\R$
defined by $t_i\mapsto x_i$ defines an archimedean ordering of $F$ where the 
$p_j$ are positive, namely an archimedean point in the given open set.

We then have that $G_{\DR}^N(k_\xi)\onto\Gamma_{k_N,{\DR}}(k_\xi)=\Z/2$ for 
each finite set $N$ of one-motives, and for all $\xi\in\Sper k$. Consequently 
the kernel of the map $H^1(k_\xi,(G_{\DR}^N)^0)\to H^1(k_\xi,G_{\DR}^N)$ is 
trivial for all $\xi$. For the dense set of $\xi\in\Sper k$ given in the 
theorem, $\eta_{N,\xi}$ is trivial in $H^1(k_\xi,G_{\DR}^N)$. Then for these 
$\xi$ we have that $\eta_{N,\xi}^0=\loc_\xi\eta_N^0$ is trivial in 
$H^1(k_\xi,(G_{\DR}^N)^0)$. 

Using the local-global principle of [Sch], Theorem 4.1, which asserts that 
for a connected linear algebraic group $G^N$ over a perfect field $k$ with 
$\vcd(k)\le 1$ the map $H^1(k,G^N)\to\prod_\xi H^1(k_\xi,G^N)$ is injective 
where the product ranges over any dense subset of orderings $\xi$ in 
$\Sper k$, we conclude that $\eta_N^0$ is $1$ for all finite sets $N$ of 
one-motives. Hence $\eta^0=\liminv\eta_N^0$ is trivial, so is its image 
$\eta$, and $P'$ and $P_\s$ define the same class in $H^1(k,G_\s)$.
\hfill\bx
\enddemo

The following lemma is used in the proof of the proposition below.

\proclaim{Lemma} Let $K_\xi$ be a real closed field containing $k_\xi$.
Then the group of $K_\xi$-points of $\Gamma_{k,{\DR}}=
{}_{z'}\Gamma_k$ is isomorphic to $\Gal(\ov k/k_\xi)$.
\endproclaim

\demo{Proof} We have 
$\Gamma_{k,\DR}(K_\xi)=\Gamma_{k,\DR}(K)^{\Gal(K/K_\xi)}$,
where $K=K_\xi(\sqrt{-1})$, and $\Gamma_{k,\DR}(K)=\Gamma_{k,\DR}(\ov k)$.
Moreover, the restriction to $\ov k$ of the non trivial element of 
$\Gal(K/K_\xi)$ is the non trivial element of $\Gal(\ov k/k_\xi)$.
The group $\Gamma_{k,\DR}$ is the profinite group scheme
attached to the identity cocycle $z'(\tau)=\tau$ in $Z^1(k,\Gamma_k)$
(this is called the Artin cocycle, see [W]). Thus $\tau\in\Gal(\ov k/k)$
acts on $\gamma\in\Gamma_{k,\DR}(\ov k)=\Gal(\ov k/k)$ by $\tau_{\DR}(\gamma)
=\tau\gamma\tau^{-1}$. In particular $c_\xi\in\Gal(\ov k/k_\xi)$ acts on
$\gamma\in\Gamma_{k,\DR}(\ov k)$ by $c_{\xi,\DR}(\gamma)
=c_\xi\gamma c_\xi^{-1}$. Hence $\Gamma_{k,\DR}(k_\xi)=\{\gamma\in
\Gal(\ov k/k); c_\xi\gamma c_\xi^{-1}=\gamma\}$. It remains to determine the 
centralizer of $c_\xi\in\Inv(\Gal(\ov k/k))$ in $\Gal(\ov k/k)$. We claim it 
is $\{1,\,c_\xi\}$. The field $k_\xi=\ov k^{c_\xi}$ of fixed points of $c_\xi$
in $\ov k$ is a real closure of $k$ whose ordering induces $\xi$ on $k$. 
If $\gamma\in\Gal(\ov k/k)$ commutes with $c_\xi$ then it maps $k_\xi$ to
itself. But the only automorphism of $k_\xi$ over $k$ is the identity (by the
Artin-Schreier theorem; see, e.g., [Sc], Ch. 3, Theorem 2.1). Hence 
$\gamma\in\Gal(\ov k/k_\xi)=\{1,\,c_\xi\}$.
\hfill\bx
\enddemo

The following proposition is used in the proof of the Theorem above. 

\proclaim{Proposition} The map $G_{\DR}(k_\xi)\to\Gamma_{k,\DR}(k_\xi)$
is surjective for every archimedean ordering $\xi$ in $\Sper k$.
\endproclaim

\demo{Proof} The lemma implies that $\Gamma_{k,\DR}(k_\xi)=\Z/2=
\Gamma_{k_\xi,\DR}(k_\xi)$. Write $G_{k,\s}$ and $G_{k,\DR}$ to specify
the base field. Using the functor $\uM_k\to\uM_{k_\xi}$ which is
induced from $M\mapsto M\times_kk_\xi$ (incidentally, it is almost surjective
(by which we mean that each object of $\uM_{k_\xi}$ is a subquotient of an
object in the image of $\uM_k$), by the proof of [DM], 6.23 (a)), we have a 
$k_\xi$-homomorphism $G_{k_\xi,\DR}\to G_{k,\DR}$ (in fact an injection, by 
[Sa], II, 4.3.2 g) ii), or [DM], Proposition 2.21 (b)) of the motivic Galois 
groups for the de Rham fiber functor.
Hence it suffices to prove the proposition only for a real closed $k$.
Since $\xi$ is archimedean, $k$ embeds in $\R$, and it suffices to prove
the proposition for $k=\R$. Thus we assume from now on that $k$ is $\R$, 
and write $G_{\DR}$ for $G_{\R,\DR}$.

Recall that the functors $\uM^0_{\R}\to\uM_{\R}$ and $\uM_{\R}\to\uM_{\C}$, 
and the fiber functor $H^{\#}_\s$, define the exact sequence $1\to G^0_\s\to
G_\s\to\Gamma_{\R}\to 1$ of affine group schemes over $\Q$ (for the ``pure''
case, which implies at once the ``mixed'' case, see [DM], Proposition
6.23(a,b)). Using the functors $\uM^0_{\R}\otimes \R\to\uM_{\R}\otimes \R$ 
and $\uM_{\R}\otimes \R\to\uM_{\C}\otimes \R$, and the fiber functor 
$H^{\#}_\s\otimes \R$, the groups become groups over $\R$ (note that [DM], 
Remark 3.12, applies with any -- not necessarily finite -- field extension 
$k'/k$). But we do not change the notations. 

For any subfield $K$ of $\R$, a $K$-{\it Hodge structure} (``over $\C$'') is
a pair $(V,(V^{p,q}))$ consisting of a finite dimensional vector 
space $V$ over $K$, and a direct sum decomposition $V\otimes_K\C=\oplus 
V^{p,q}$ with $\tau_\infty(V^{p,q})=V^{q,p}$; $\tau_\infty\not=1$ in
$\Gal(\C/\R)$. A $K$-{\it Hodge structure over} $\R$ is a triple 
$(V,(V^{p,q}),F_\infty)$ where the new ingredient is an involutive 
endomorphism $F_\infty$ of $V$ whose extension to $V\otimes_K\C$ satisfies 
$F_\infty(V^{p,q})=V^{q,p}$. With the natural definition of tensor products 
and morphisms, these make neutral Tannakian categories $\uHod_K$ ($K$-Hodge 
structures) and $\uHod_K^+$ ($K$-Hodge structures over $\R$) over $K$ (for 
the forgetful fiber functor $\omega_K:(V,\dots)\to V$). 

A $K$-{\it mixed Hodge structure} (``over $\C$'') is a triple $(V,W_{\bullet},
F^{\bullet})$, where $V$ is a finite dimensional $K$-vector space with a 
finite increasing (weight) filtration $W_{\bullet}$ and a finite decreasing 
(Hodge) filtration $F^{\bullet}$ on $V\otimes_K\C$, such that $F^{\bullet}$ 
induces a $K$-Hodge structure of weight $n$ on the graded piece 
$\Gr_n^WV=W_nV/W_{n-1}V$ for each $n$. A $K$-{\it mixed Hodge structure over} 
$\R$ is a $K$-mixed Hodge structure $(V,W_{\bullet},F^{\bullet})$ with a 
$W_{\bullet}$ preserving involutive automorphism $F_\infty$ of $V$ such that 
$F_\infty((\Gr_n^WV\otimes_K\C)^{p,q})=(\Gr_n^WV\otimes_K\C)^{q,p}$. With the
natural definition of $\otimes$ and morphisms, these make the Tannakian
categories $\MHS_K$ and $\MHS^+_K$.

The main Theorem 2.11 of [D2] asserts that for an algebraically closed 
subfield $\frak K$ of $\C$, the functor $H_\s:\uM_{\frak K}^{\red}\to\uHod_\Q$
is fully faithful. It is extended in [D-III], 10.1.3, to assert that the 
functor $H_\s:M\mapsto H_\s(M)=\s M(\C)_B$ defines an equivalence between the 
category of isogeny classes of one-motives over $\frak K$ and the category of 
($\Q$-)mixed Hodge structures of type $\{(0,0),\,(0,-1),\,(-1,0),\,(-1,-1)\}$
whose graded quotient $\Gr_{-1}$ is polarizable. 
A morphism of one-motives is a morphism $(\alpha,\beta)$ of complexes 
$[X\to E]\to[X'\to E']$. It is an {\it isogeny} if both $\alpha$ and $\beta$
are isogenies, i. e. have finite kernels and cokernels.
The functor $H_{\s}$ extends to a 
faithful functor from the Tannakian
category $\uM_{\C}$ to the Tannakian category $\MHS_{\Q}$ (in this context
we note Theorem 2.2.5 of [Br], which asserts that a Hodge cycle on a 
one-motive -- and in particular a power thereof -- is absolute), and from
$\uM_{\R}$ to $\MHS_{\Q}^+$: $\tau_\infty\in\Gal(\C/\R)$ induces an involution
of $\s M(\C)$, hence an involution $F_\infty=H_\s(\tau_\infty)$ on $H_\s(M)$.
The restriction of $H_{\s}$ to $\uM^0_{\R}$ is an equivalence with the 
category $\Rep_{\Q}\Gamma_{\R}$ of representations of $\Gamma_{\R}$ over $\Q$.

The fiber functor $H^{\#}_{\s}\otimes\R$ on $\uM_{\R}\otimes\R$ factorizes 
through the forgetful functor 
$\omega_{\R}:\MHS^+_{\R}\to\Rep_{\R}\Gamma_{\R}$. The 
restriction of $\omega_{\R}$ to $\MHS_{\R}$ is the forgetful functor into
the category $\Vec_{\R}$ of vector spaces over $\R$. The restriction of 
$H^{\#}_{\s}\otimes \R$ to $\uM^0_{\R}\otimes \R$ is an equivalence of 
categories with $\Rep_{\R}\Gamma_{\R}$. 

But we are concerned with the fiber functor $H^{\#}_{\DR}\otimes \R$ 
on $\uM_{\R} \otimes \R$ and the exact sequence 
$1\to G^0_{\DR}\to G_{\DR}\to \Gamma_{\R,\DR}\to 1$
of real groups associated with the almost surjective functor 
$\uM_{\R}\otimes\R\to\uM_{\C}\otimes\R$ and the fully faithful functor
$\uM^0_{\R}\otimes\R\to\uM_{\R}\otimes\R$. To show the surjectivity of the
map $G_{\DR}(\R)\to \Gamma_{\R,\DR}(\R)$ of groups of real points, it suffices
to show that the reductive part $G^{\red}_{\DR}(\R)$ surjects on 
$\Gamma_{\R,\DR}(\R)$. For this, note that the functor $H^{\#}_{\DR}\otimes\R$
on $\uM^{\red}_{\R}\otimes\R$ factorizes via $\uM^{\red}_{\R}\otimes\R\to
\Hod^+_{\R}$ and a functor $\omega_{\DR,\R}:\Hod^+_{\R}\to\Vec_\R$
described below. This follows from the fact that for the realizations
of a motive one has $c_{\DR}=F_{\infty}\circ\br$, where $c_{\DR}$ and $\br$
are respectively the deRham and the Betti complex conjugations.
Defining $\Bbb S^+_{\DR}=\Aut^\otimes(\omega_{\DR,\R}|
\Hod^+_{\R})$ (and $\Bbb S_{\DR}=\Aut^\otimes(\omega_{\DR,\R}|\Hod_{\R})$),
we get the vertical arrow in the commutative square
$$\matrix \Bbb S^+_{\DR}&\to&\Gamma_{\R,\DR}\\ \downarrow&&\Vert\\
G^{\red}_{\DR}&\to&\Gamma_{\R,\DR}\,.\endmatrix$$
The horizontal arrows result from the fully faithful functors 
$\uM^0_{\R}\otimes\R\to\uM^{\red}_{\R}\otimes\R$ and $\uM^0_{\R}\otimes\R\to
\Hod^+_{\R}$. Consequently it suffices to show that $\Bbb S^+_{\DR}(\R)\onto
\Gamma_{\R,\DR}(\R)$.

Analogously we have the functor $\omega_{\R}$ on $\Hod^+_{\R}$, the
real groups $\Bbb S^+=\Aut^\otimes(\omega_{\R}|\Hod^+_{\R})$ and 
$\Bbb S=\Aut^\otimes(\omega_{\R}|\Hod_{\R})$, and the exact sequence
$1\to \Bbb S\to \Bbb S^+\to \Gamma_{\R}\to 1$. The motivic Galois
group $\Bbb S$ of $\Hod_{\R}$ and the functor $\omega_{\R}$ is well 
known ([DM], Example 2.31). The group $\Bbb S$ is the connected
$\R$-group $\Res_{\C/\R}\G_m$ obtained from the multiplicative group 
$\G_m$ on restricting scalars from $\C$ to $\R$. Thus $\Bbb S(\C)=
\C^\times\times\C^\times$, and the non-trivial element of $\Gal(\C/\R)$ 
acts on $\Bbb S(\C)$ by $(a,b)\mapsto (\ov b,\ov a)$, so $\Bbb S(\R)
=\C^\times$. Indeed, a representation $\rho:\Bbb S\to\Aut(V)$ defines 
$V^{p,q}$ to be the $v\in V\otimes_\R\C$ with $\rho(z)(v)=z^{-p}\ov z^{-q}v$ 
for all $z\in\C^\times$. The motivic Galois group of the subcategory 
$\uHod_\R^0$ of the $V$ in $\uHod_\R^+$ with $V^{p,q}=\{0\}$ unless $p=q=0$ 
is the constant group scheme $\Gamma_\R$ over $\R$ associated to the group 
$\Gal(\C/\R)$. The motivic Galois group of $\uHod_\R^+$ (and $\omega_\R$) is 
an extension $\Bbb S^+$ of $\Gamma_\R$ by $\Bbb S$.
Indeed, a triple $(V,(V^{p,q}),F_\infty)$ is associated with the extension of
$\rho$ from $\Bbb S$ to $\Bbb S^+$ by $\rho(1\times\br)=F_\infty$ (``bar'' 
signifies complex conjugation).
The exact sequence $1\to \Bbb S\to \Bbb S^+\to\Gamma_\R\to 1$ is defined by 
the fully faithful functor $\uHod_\R^0\to\uHod_\R^+$ and the essentially 
surjective ``forget $F_\infty$'' functor $\Hod_\R^+\to\uHod_\R$.
Note that the sequence is split, and $\Bbb S^+=\Bbb S\ltimes\Gamma_{\R}$.
A splitting is given by the essentially surjective functor 
$\Hod_\R^+\to\uHod_\R^0$, ``forget the Hodge structure'', and $\Gamma_{\R}$
acts on $\Bbb S$ via the Galois action.

Since $H^1(\R,\Bbb S)=1$, the sequence 
$1\to \Bbb S(\R)\to \Bbb S^+(\R)\to \Gamma_{\R}(\R)\to 1$ is exact. Since 
the group $\Bbb S_{\DR}$ is $\G_m^2$ (see the following paragraph), by Hilbert
Theorem 90 we have $H^1(\R,\Bbb S_{\DR})=1$. Hence $\Bbb S^+_{\DR}(\R)\to 
\Gamma_{\R,\DR}(\R)$, which is just 
$H^0(\R,S^+_{\DR})\to H^0(\R,\Gamma_{\R,\DR})$, is onto.
This completes the proof of the proposition.

Note that the structure of the entire group
$\Aut^\otimes(\omega_{\DR,\R}|\MHS_{\R})$
is computed in [D3], Construction 1.6 and Proposition 2.1, since 
$\omega_{\DR,\R}$ is the functor $\Gr^W$ of [D3]. But by the Levi 
decomposition it suffices for us to work only with its reductive part.
Thus we note that $\Bbb S^+_{\DR}$ is known to be 
$(\G_m\times\G_m)\rtimes\Z/2$. Indeed,
the category $\uHod_\R^+$ is equivalent to the category $\uHod_\R^\ast$ of 
triples $(W,(W^{p,q}),F)$, where $W$ is a finite dimensional real vector space
with decomposition $W=\oplus W^{p,q}$ into real subspaces, and $F$ is an 
involutive endomorphism of $W$ over $\R$ with $F(W^{p,q})=W^{q,p}$. In fact, 
$\uHod_\R^+\to\uHod_\R^\ast$ is given by $W^{p,q}=$ fixed points of 
$F_\infty\circ\br$ in $V^{p,q}$, $F=F_\infty|W$, $W^{p,q}=W\cap V^{p,q}$, 
and $\uHod_\R^\ast\to\uHod_\R^+$ by: $V=$ fixed points of $F\circ\br$ in 
$W\otimes\C$, $V^{p,q}=V\cap(W^{p,q}\otimes\C)$, $F_\infty=F|V$. 
The fiber functor $H^{\#}_{\DR}\otimes\R$ on $\uM_{\R}\otimes\R$ factorizes 
through the fiber functor $\omega_{\DR}$ on $\uHod_\R^+$, which is $V\mapsto
W$, or $W\mapsto W$ on $\uHod_\R^\ast$. The group of automorphisms of
$\omega_{\DR}$ on $\uHod_\R^\ast$ is 
$\Bbb S^+_{\DR}=(\G_m\times\G_m)\rtimes\Z/2$, the product of the finite group 
scheme $\Gamma_{\R,\DR}=\Z/2$ by $\Bbb S_{\DR}=\G_m\times\G_m$, the groups of 
automorphisms of 
the functor $\omega_{\DR}$ on the categories $\uHod_\R^0$ and $\uHod_\R$.
\hfill\bx
\enddemo

\remark{Remark 1} Proposition 6.22(b) of [DM] is wrong ($\unM_k\to\unM_{k'}$ 
there is fully faithful but not essentially surjective), but this is of no 
consequence for the theory. For a corrected statement and a counter example
see [S2], \S6. The connectedness assertion in Proposition 6.22(a)
(and consequently 6.23) of [DM] -- which is a consequence of the standard 
conjectures -- is out of reach of current technology (in Deligne's opinion) in
the context of the whole category of (even only pure) motives. In particular, 
(6.1) of [S2] should be (6.1?), and similarly for [J], Theorem 4.7, p. 50. 
The proof of [DM], 6.22(a), implicitly assumes that Hodge cycles are absolute.
It works in our setting (of motives of abelian varieties, and one-motives) 
since Hodge cycles on abelian varieties are absolute, by [D2], Theorem 2.11. 
Thus we use [DM], 6.22(a) and 6.23, replacing $\unM_k$, $\unM_{\ov k}$ by 
$\uM_k^{\red}$, $\uM_{\ov k}^{\red}$ 
in [DM], p. 213, l. -7 to p. 216, l. -9; in particular the group $G(\s)$ 
of [DM], p. 213, l. -6 (denoted $G_\s$ here) should be 
$\uAut^\otimes(H^{\#}_\s|\uM_k^{\red})$, and in the proof of 
[DM], 6.22(a), $X$ should be in $\uM_k^{\red}$ (to use (I 3.4)). 

Yet the full Galois group $G(\s)$ of [DM], 6.22(a) 
($=\uAut^\otimes(H^{\#}_\s|\unM_{\ov k})$) is pro-reductive (as asserted in 
[DM], 6.22(a)) -- meaning that its connected component $G^0$ is the projective
limit of connected reductive groups -- by [DM], Remark 2.28 (``$G^0$ is 
pro-reductive iff $\uRep_{\Q}G(\s)$ is semi-simple'') and [DM], 
Proposition 6.5 (``$\unM_{\ov k}=\uRep_{\Q}G(\s)$ is semi-simple'').
\endremark

In an attempt to clarify the proof of [DM], 6.22(a), note that it uses the 
following well-known assertion. Only the special case of pure Hodge
structures is used in [DM], and this suffices for our purposes too,
since an algebraic group is connected if its (Levi) reductive component is. 
As in [DM], Proof of Proposition 2.8, let $C_H$ be the full (Tannakian) 
subcategory of the category $\Hod$ of $\Q$-Hodge structures generated by 
$\Q(1)$ and an object $H$. The objects of $C_H$ are by definition
the subquotients of sums of $T=H^{\otimes m_1}\otimes (H^{\sssize\vee})
^{\otimes m_2}\otimes\Q(1)^{\otimes m_3}$, and $a\in\G_m$ acts on 
$\Q(1)^{\otimes m}$ by multiplication by $a^{-m}$. 
Let $\omega$ be the fiber (forgetful) functor to the category of vector 
spaces over $\Q$. 
Suppose that $H$ is a polarizable Hodge structure. Then $C_H$ is semi-simple.
Write $G'$ for the subgroup $GL(H)\times\G_m$ over $\Q$ which fixes all 
$(0,0)$-vectors $t$ in every object $T$ of $C_H$.

\proclaim{Assertion} The group 
$G=\operatorname{Aut}^\otimes(\omega\vert C_H)$ 
is isomorphic to the group $G'$.
\endproclaim

\demo{Proof} A morphism $g=(g_X:\Phi(X)\to \Phi'(X))$ of functors $\Phi$, 
$\Phi'$ on a category satisfies $\Phi'(f)g_X=g_Y\Phi(f)$ for every morphism 
$f:X\to Y$. In $C_H$, an endomorphism of the fiber functor 
$\omega$ is an element $g$ of $GL(H)\times\G_m$ which -- extended to 
$H_{\C}=H\otimes\C$ -- commutes with $\omega(f)$, thus $g\omega(f)=
\omega(f)g$, for every morphism $f:V\to U$ in $\Hod$, namely with all 
linear maps $f:V\to U$ with $f(V^{p,q})\subset U^{p,q}$. Thus for each $V$, 
$g$ commutes with $\Hom_{\Hod}(\Q(0),V)=V^{0,0}$, namely it fixes $V^{0,0}$,
so $g\in G'$.

Conversely, if $g\in G'$ then for any $V,\,U\in C_H$, $g$ fixes
$(V^{\sssize\vee}\otimes U)^{0,0}=\Hom(V,U)^{0,0}$, thus $g:H\to H$
commutes with every morphism $f:V\to U$ in $\Hod$, so $g\in G$.
\hfill\bx
\enddemo

Now the problem in the proof of 6.22(a) in [DM] is that for $X$ in the 
Tannakian category $\unM_{\frak K}$ of motives of absolute Hodge cycles, the 
full subcategory $C_X$ of $\unM_{\frak K}$ embeds via $H_\s$ in the Tannakian
category $\Hod$ of $\Q$-Hodge structures, but it is not a full subcategory
unless each $\s$-Hodge cycle on $X$ is absolute. If $C_X$ is a full 
subcategory of $\Hod$ (via $H_\s$, namely each $\s$-Hodge cycle is absolute), 
then $G_X=\Aut^\otimes(H_\s\vert C_X)$ of [DM], 6.22(a), becomes the group $G$
of the Assertion above, and it can be compared with $G'$, the connected
group which features in the second half of [DM], proof of 6.22(a) (and (I 3.4)
there). In general, the group $G_X$ consists of those automorphisms of the 
vector space $H_\s(X)$ which commute with each automorphism of the absolute
Hodge structure $H(X)$. Not every automorphism $f_\s$ of the Hodge 
structure $H_\s(X)$ extends to an automorphism $(f_{\DR},f_\ell,f_\tau)$ 
of absolute Hodge structures, so the group $G_X$ -- being the commutator of
absolute Hodge morphisms -- may be larger than the commutator $G$ of the
larger family of $\s$-Hodge morphisms. The two groups are equal (and the 
{\it a-priori} possibly bigger $G_X$ is connected) for abelian varieties $X$, 
for which Hodge cycles are absolute.

\remark{Remark 2} An extension $E$ of an abelian variety by a torus $T$ is 
commutative: (a) $T$ is central: the action by inner automorphism of $A=E/T$ 
on $T$ is trivial, because it amounts to an action on the character group,
which is discrete; (b) the commutator $E\times E\to E$ has image in $T=\ker[E
\to A]$, and it factors via $A\times A=E/T\times E/T\to T$ by (a); it is 
trivial since the image is proper and reduced in the affine $T$.
\endremark

\bigskip\bigskip
\mathsurround 2pt
\def\ref#1#2{\noindent\hangindent 4em\hangafter1\hbox to 4em{\hfil#1\quad}#2}

\n {\bf References}
\medskip
\ref{[BP]}{E. Bayer-Fluckiger, R. Parimala, Classical groups and the Hasse 
principle, {\it Ann. of Math.} 147 (1998), 651-693.}

\ref{[BB]}{D. Blasius, M. Borovoi, On period torsors, {\it Automorphic forms, 
automorphic representations, and arithmetic}, Proc Symp. Pure Math. 66, Part 
1, AMS, Providence, RI, 1999.}

\ref{[BCR]}{J. Bochnak, M. Coste, M.-F. Roy, {\it G\'eom\'etrie alg\'ebrique 
r\'eelle}, Ergebnisse der Math. 12, Springer-Verlag (1987).}

\ref{[Br]}{J.-L. Brylinski, ``1-Motifs'' et formes automorphes, in {\it
Journ\'ees Automorphes}, Publ. Math. Univ. Paris VII 15 (1983), 43-106.}

\ref{[D-II]}{P. Deligne, Th\'eorie de Hodge, II, {\it Publ. Math. IHES} 40 
(1971), 5-58.}

\ref{[D-III]}{P. Deligne, Th\'eorie de Hodge, III, {\it Publ. Math. IHES} 44 
(1975), 5-77.}

\ref{[D1]}{P. Deligne, Vari\'et\'es de Shimura: interpr\'etation modulaire, et
techniques de construction de mod\`eles canoniques, in {\it Automorphic Forms,
Representations and $L$-functions}, Proc. Sympos. Pure Math. 33 II (1979), 
247-290.}

\ref{[D2]}{P. Deligne, (Notes by J. Milne), Hodge cycles on abelian varieties,
in {\it Hodge Cycles, Motives, and Shimura Varieties}, Lecture Notes in 
Mathematics 900, Springer-Verlag (1982), 9-100.}

\ref{[D3]}{P. Deligne, Structures de Hodge mixtes r\'eelles, in {\it Motives},
Proc. Sympos. Pure Math. 55 (1994), 509-514.}

\ref{[DM]}{P. Deligne, J. Milne, Tannakian categories, in {\it Hodge Cycles,
Motives, and Shimura Varieties}, Lecture Notes in Mathematics 900,
Springer-Verlag (1982), 101-228.}

\ref{[J]}{U. Jannsen,  Mixed motives and algebraic K-theory, Lecture Notes 
in Mathematics 1400, Springer-Verlag (1990).}


\ref{[M]}{J. Milne, Canonical models of (mixed) Shimura varieties and
automorphic vector bundles, in  {\it Automorphic Forms, Shimura Varieties,
and $L$-functions} I (1990), 283-414.}

\ref{[Sa]}{N. Saavedra Rivano, {\it Cat\'egories Tannakiennes}, Lecture Notes 
in Mathematics 265, Springer-Verlag (1972).}

\ref{[Sc]}{W. Scharlau, {\it Quadratic and Hermitian Forms}, Grundlehren 270,
Springer-Verlag (1985).}

\ref{[Sch]}{C. Scheiderer, Hasse principles and approximation theorems for
homogeneous spaces over fields of virtual cohomological dimension one, {\it
Invent. Math.} 125 (1996), 307-365.}

\ref{[S1]}{J.-P. Serre, {\it Cohomologie Galoisienne}, Cinqui\`eme \'edition,
Lecture Notes in Mathematics 5, Springer-Verlag (1994).}

\ref{[S2]}{J.-P. Serre, Propri\'et\'es conjecturales des groupes de Galois
motiviques et des\hb repr\'esentations $\ell$-adiques, in {\it Motives}, Proc.
Sympos. Pure Math. 55 I (1994), 377-400.}


\ref{[W]}{J.-P. Wintenberger, Torseurs pour les motifs et pour les
repr\'esentations $p$-adiques potentiellement de type CM, {\it Math. Ann.}
288 (1990), 1-8.}
\enddocument

\bye